\documentclass[12pt]{amsart}
\usepackage{graphics}
\usepackage{amssymb}
\usepackage{amsthm}
\usepackage{pictex}

\begin{document}

\newtheorem{lemma}{Lemma}[section]
\newtheorem{prop}[lemma]{Proposition}
\newtheorem{cor}[lemma]{Corollary}
\newtheorem{thm}{Theorem}
\newtheorem*{thm*}{Theorem}

\newtheorem{rem}[lemma]{Remark}
\newtheorem{rems}[lemma]{Remarks}
\newtheorem{defi}[lemma]{Definition}
\newtheorem{ex}{Example}
\newtheorem{convention}[lemma]{Convention}

\newcommand{\ra}{\longrightarrow}
\newcommand{\cstar}{$C^\ast$-algebra }
\newcommand{\cstars}{$C^\ast$-algebras }
\newcommand{\cH}{\mathcal H}
\newcommand{\oh}{\overline{h}}
\newcommand{\ok}{\overline{k}}
\newcommand{\St}{\hbox{st}}
\newcommand{\lk}{\hbox{lk}}

\title[Ozawa kernels]{
Exactness of free and amenable groups by the construction of Ozawa
kernels }
\author{Sarah Campbell $^\dagger$}
\date{\today}
\thanks{$^\dagger$ Supported by EPSRC}
\begin{abstract}  Using properties of their Cayley graphs, specific
examples of Ozawa kernels are constructed for both free and
amenable groups, thus showing that these groups satisfy Property
$O$. It is deduced both that these groups are exact and satisfy
Yu's Property $A$.
\end{abstract}

\maketitle
\section{Introduction}

The property of exactness was first introduced as an analytic
property of $C^\ast$-algebras. We say that a $C^\ast$-algebra $C$
is exact if the operation of taking the cross product with $C$
preserves exactness of short exact sequences.  In other words, for
any short exact sequence $$0\rightarrow J\rightarrow B\rightarrow
B/J\rightarrow0$$ the sequence $$0\rightarrow J\otimes_{min}
C\rightarrow B \otimes_{min} C\rightarrow B/J \otimes_{min} C
\rightarrow0$$ is also exact.

We say that a group $\Gamma$ is exact if the operation of taking
the reduced crossed product with $\Gamma$ preserves exactness of
short exact sequences of $\Gamma$-$C^*$-algebras. In other words,
$\Gamma$ is exact if and only if for every exact sequence of
$\Gamma$-$C^*$-algebras
\[0 \ra B \ra C \ra D\ra 0 \] the sequence \[ 0 \ra C^*_r(\Gamma,
B) \ra C^*_r(\Gamma, C) \ra C^*_r(\Gamma , D)\ra 0 \] of crossed
product algebras is exact. Following work by Kirchberg and
Wassermann, \cite{KW}, a discrete group is said to be exact if and
only if its reduced \cstar is exact. Exactness is linked to
conjectures such as the Novikov Conjecture and the Baum Connes
conjecture and exact groups satisfy the Coarse Baum Connes
Conjecture. Examples of exact groups include hyperbolic groups,
groups with finite asymptotic dimension and groups acting on
finite dimensional CAT(0) cube complexes \cite{CN}.

More recent work has shown that when we regard the group as a
metric space in the word metric, the property of exactness can be
defined by more geometric means. In particular, Yu showed that
exactness of a group is equivalent to Property $A$, a measure
theoretic property, and implies the Uniform Embedding Property
\cite{Yu}. In \cite{Ozawa}, Ozawa introduced the following
property which we will call Property $O$ and proved that for a
discrete group it is equivalent to exactness of its reduced
$C^*$-algebra:
\begin{defi} A discrete group $G$ is said to have Ozawa's Property O
if for any finite subset $E\subset G$ and any $\epsilon >0$, there
are a finite subset $F\subset G$ and $u$: $G\times G
\rightarrow\mathbb{R}$ such that
\begin{enumerate}
\item $u(x,y)$ is a positive definite kernel
\item $u(x,y)\neq 0$ only if $x^{-1}y\in F$
\item $|1-u(x,y)|<\epsilon$ if $x^{-1}y\in E$
\end{enumerate}\end{defi}

We will call functions satisfying the conditions of Property $O$
Ozawa kernels.

This property has been used by Guentner and Kaminker to prove
their theorem relating asymptotic compression and exactness
\cite{GK}. No explicit examples of Ozawa kernels can be found in
the literature. The aim of this note is to construct explicit
Ozawa kernels for two classical cases of exact groups, amenable
groups and free groups by using geometric properties of their
Cayley graphs. In both cases the functions can be viewed as
weighted mass functions.  This shows that amenable and free groups
satisfy Property O. This is sufficient for us to deduce both that
their reduced \cstar is exact and that they satisfy Yu's property
$A$.
\section{Free groups}

\begin{thm}Free groups admit an Ozawa kernel and thus satisfy
Property $O$.\end{thm}

\begin{proof}

This construction is based on the proof that trees have Property A
\cite{Janus}.

Let $T$ be the Cayley graph of a free group (a tree) and $V$ its
set of vertices. Let $\gamma_0: \mathbb{R}\rightarrow T$ be a
geodesic ray in $T$. Let $\gamma_v$ be the unique geodesic ray
issuing from $v\in V$ and intersecting $\gamma_0$ along a geodesic
ray. Let $\gamma_v^n$ be the initial segment of $\gamma_v$ of
length $n$.

For any $x,y \in V$, we define our function $u_n(x,y)$  to be the
size of the overlap of the n-length rays $\gamma_x^n$ and
$\gamma_y^n$, scaled by $n+1$.
\[u_n(x,y)=\frac{|\gamma_x^n\cap\gamma_y^n|}{n+1}\] We will now show that this family of functions can be used to
define an Ozawa kernel.

\begin{lemma}For any $n$, $u_n$ is a positive definite kernel.\end{lemma}

\begin{proof}
Define $f_v(x)$ and $\chi_{v,n}(x)$ as follows:
\begin{displaymath}
f_v(x)= \left\{
\begin{array}{ll} 1 &
\textrm{if $v$ separates $x$ from the end of $\gamma_0$}
\\ 0& \textrm{otherwise}\\
\end{array} \right.
\end{displaymath}
And
\begin{displaymath}
\chi_{v,n}(x)= \left\{
\begin{array}{ll} 1& \textrm{if $x$ is contained in the ball of radius $n$ around
$v$}\\ 0 & \textrm{otherwise}
\end{array} \right.
\end{displaymath}

Then we have \[u_n(x,y)=\frac{1}{n+1}\sum_{v\in V}f_v(x) f_v(y)
\chi_{v,n}(x) \chi_{v,n}(y)\]

The only vertices contributing to this sum are whose which are
within distance $n$ of both $x$ and $y$ and which separate both
$x$ and $y$ from the end of $\gamma_0$.  This is precisely the
size of the intersection of the $n$-length rays $\gamma_x^n$ and
$\gamma_y^n$.

We can now rearrange $u_n$ to show that this is a positive
definite kernel.

\begin{eqnarray*} \sum_{i,j}^{k} \lambda_i\lambda_j
u_n(x_i,x_j)&=& \sum_{i,j}^{k}\lambda_i\lambda_j \sum_{v\in
V}f_v(x_i)f_v(x_j) \chi_{v,n}(x_i)\chi_{v_{B_n}}(x_j)\\ &=&
\sum_{v\in V} \left(\sum_i^k \lambda_i
f_v(x_i){\chi_{v,n}}(x_i)\sum_j^k \lambda_j
f_v(x_j){\chi_{v,n}}(x_j)\right)\\ &=& \sum_{v\in V}\left(\sum_i^k
\lambda_i f_v(x_i){\chi_{v,n}}(x_i)\right)^2 \geq 0\\
\end{eqnarray*}
\end{proof}

\begin{lemma}For each $n$ there exists a finite set $F$ such that
$u_n(x,y)\neq 0$ only if $x^{-1}y\in F$.\end{lemma}

\begin{proof}Let $F$ be the ball of radius $n$ around the origin. If
$d(x,y)>2n$, ie $x^{-1}y\notin F$ there is no overlap between the
$n$-length rays $\gamma_x^n$ and $\gamma_y^n$ and so
$u_n(x,y)=0$.\end{proof}

\begin{lemma}Given any finite subset $E$ and $\epsilon>0$ there exists
$N$ such that  $|1-u_N(x,y)|<\epsilon$ if $x^{-1}y\in
E$.\end{lemma}

\begin{proof}Since $E$ is a finite subset, there exists some number $m$
such that if $x^{-1}y \in E$ then $d(x,y)<m$.

Now if $d(x,y)<m$, then the minimum size of the overlap of
$\gamma_x^n$ and $\gamma_y^n$ is $n-m$ whilst the maximum is
$n+1$.

And so we have \[\frac{n-m}{n+1}\leq u_n(x,y) \leq
\frac{n+1}{n+1}\]

Hence $\lim_{n\rightarrow \infty}u_n(x,y)=1$

And so $\forall \epsilon>0$, $\exists N$ such that $\forall
x^{-1}y \in E$ we have as required
$|1-u_N(x,y)|<\epsilon$.\end{proof} Hence $u_N$ is an Ozawa type
kernel and free groups satisfy Property $O$.
\end{proof}

The following corollary immediately follows:
\begin{cor}Free groups are exact and satisfy Yu's property $A$.\end{cor}

\pagebreak

\section{Amenable groups}

We will first review the definition of an amenable group via
Folner's condition:
\begin{defi}If a group $G$ is amenable, then there exists a Folner
sequence ${G_n}$ of finite subsets of $G$ such that $\forall g \in
G$, \[\lim_{n\rightarrow \infty}\frac{|gG_n\vartriangle
G_n|}{|G_n|}=0\]\end{defi}

This can be rewritten as follows:
\begin{eqnarray*}
\frac{\left|gG_n\vartriangle
G_n\right|}{\left|G_n\right|}&=&\frac{|gG_n\cup G_n-gG_n\cap
G_n|}{|G_n|}\\ &=&\frac{|gG_n\cup G_n|}{|G_n|}-\frac{|gG_n\cap
G_n|}{|G_n|}
\end{eqnarray*}

The maximum possible value of $|gG_n\cup G_n|$ is $2|G_n|$ and the
minimum is $|G_n|$.  Similarly, the maximum possible value of of
$|gG_n\cap G_n|$ is $|G_n|$ and the minimum is 0.

So $\frac{|gG_n\cup G_n|}{|G_n|}$ lies between 1 and 2, while
$\frac{|gG_n\cap G_n|}{|G_n|}$ lies between 0 and 1.  Since the
difference between them tends to 0 and we are dealing with bounded
sets of real numbers, we have that $\lim_{n\rightarrow \infty}
\frac{|gG_n\cap G_n|}{|G_n|}=1$.

We can now prove our next theorem:
\begin{thm}Amenable groups admit an Ozawa kernel and thus
satisfy Property $O$.\end{thm}

\begin{proof}
Consider the family of functions $u_n(x,y)= \frac{|xG_n \cap
yG_n|}{|G_n|}$.

\begin{lemma}For each $n$, $u_n$ is a positive definite kernel.\end{lemma}

\begin{proof}An element $g\in G$ belongs to the intersection $xG_n\cap
yG_n$ only if $g\in xG_n$ and $g\in yG_n$. This is equivalent to
$x^{-1}g\in G_n$ and $y^{-1}g\in G_n$.  Again this is equivalent
to $x^{-1}\in G_ng^{-1}$ and $y^{-1}\in G_ng^{-1}$.   And so the
condition becomes $x\in gG_n^{-1}$ and $y\in gG_n^{-1}$.

So taking $\chi(x)$ to be the characteristic function, $u_n$ can
be rewritten as \[u_n(x,y)=\frac{1}{|G_n|}\sum_{g\in
G}\chi_{gG_n^{-1}}(x)\chi_{gG_n^{-1}}(y)\] Thus:
\begin{eqnarray*}
\sum_{i,j}^{n} \lambda_i\lambda_j u(x_i,x_j)&=&
\sum_{i,j}^{n}\lambda_i\lambda_j \frac{1}{|G_n|}\sum_{g\in
G}\chi_{gG_n^{-1}}(x_i)\chi_{gG_n^{-1}}(x_j)\\ &=&
\frac{1}{|G_n|}\sum_{g\in G} \left(\sum_i^n \lambda_i
\chi_{gG_n^{-1}}(x_i)\sum_j^n
\lambda_j\chi_{gG_n^{-1}}(x_j)\right)\\
&=&\frac{1}{|G_n|}\sum_{g\in G}\left(\sum_i^n \lambda_i
\chi_{gG_n^{-1}}(x_i)\right)^2  \geq0 \\
\end{eqnarray*}
\end{proof}

\begin{lemma}For each $n$, there exists a finite set $F$ such that $u_n(x,y)\neq 0$
only if $x^{-1}y\in F$.\end{lemma}

\begin{proof}Consider the Cayley graph of $G$.  Since $G_n$ is finite it
is contained within a ball of diameter $r$. Let $F$ be the ball of
radius $r$ around the origin. If $d(x,y)>2r$, ie $x^{-1}y\notin
F$, there is no intersection between $xG_n$ and $yG_n$ and so
$u_n=0$ as required.\end{proof}

\begin{lemma}Given any finite subset $E$ and $\epsilon>0$ there exists a
function $u_N$ such that  $|1-u_N(x,y)|<\epsilon$ if $x^{-1}y\in
E$.\end{lemma}

\begin{proof}$\lim_{n\rightarrow \infty} \frac{|gG_n\cap G_n|}{|G_n|}=1$.
So for a given $g$, $\forall \epsilon>0$, there exists $R$ such
that if $n>R$, $|1-\frac{|gG_n\cap G_n|}{|G_n|}|<\epsilon$. Since
this holds for any $g \in G$, it holds in particular for
$x^{-1}y$. Since $E$ is a finite subset, there exists
$N=\max\{R|g=x^{-1}y \in E\}$ such that $\forall x^{-1}y \in E$,
$|1-\frac{|x^{-1}yG_N\cap G_N|}{|G_N|}|=|1-\frac{|xG_N\cap
yG_N|}{|G_N|}|=|1-u_N(x,y)|<\epsilon$.\end{proof}

This function $u_N$ is an Ozawa kernel and thus amenable groups
satisfy Property $O$.
\end{proof}

The following corollary immediately follows:

\begin{cor}Amenable groups are exact and satisfy Yu's property $A$.\end{cor}

\subsection{Example: groups of subexponential growth}

A good example of the above construction which clearly shows the
importance of the geometry of the Cayley graph is that of groups
of subexponential growth which we define as follows:

\begin{defi} Let $G$ be a group with generating set $A$. Let
$\beta_A(n)$ be the number of vertices in the closed ball of
radius $n$ about 1 in the Cayley graph of the group generated by
$A$.  The growth function of $G$ with respect to $A$ is
$n\rightarrow \beta_A(n)$.\end{defi}

\begin{defi} $G$ has subexponential growth if $\beta_A(n)\leq
e^{\sqrt{n}}$ for all $n\in \mathbb{N}$.\end{defi}

Examples of groups of subexponential growth include finite groups,
abelian groups and nilpotent groups. All groups of subexponential
growth are amenable and so satisfy Folners condition. In fact, it
can be shown that balls of radius $n$ in the Cayley graph of $G$
are Folner sets \cite{BHV}.

In this case, the Ozawa kernel $u(x,y)$ is simply the size of the
intersection of the balls of radius $n$ centred at $x$ and $y$,
scaled by the size of $B_n$.  We choose the radius $n$ according
to the given $\epsilon$ and finite set $E$.
 \[u_n(x,y)=\frac{|xB_n\cap yB_n|}{|B_n|}\]
 Alternatively we can regard the function as the number of balls of radius
 $n$ which contain both $x$ and $y$, scaled by the size of a ball of radius $n$. \[u_n(x,y)=\frac{|\{B_n  | B_n \textrm{ contains both $x$ and
 $y$.\}}|}{|B_n|}\]

\end{document}